\documentstyle[12pt]{article}
\newtheorem{thm} {Theorem}
\newtheorem{cor} {Corollary}
\newtheorem{prop}[thm]{Proposition}
\newtheorem{lem}[thm]{Lemma}
\newenvironment{rem}{\vspace{3mm}\noindent {\bf Remark} }{\vspace{3mm}}
\newtheorem{conj}{Conjecture }
\newcommand{\Pf}{\noindent {\it Proof}}
\newcommand{\E}{{\rm Ext}}
\newcommand{\ra}{\rightarrow}
\newcommand{\ot}{\otimes}
\newcommand{\Ra}{\Rightarrow}
\renewcommand{\O}{{\cal O}}
\newcommand{\G}{{\Gamma}}
\newcommand{\qed}{\ \vrule width2mm height3mm \vspace{3mm}}
\renewcommand{\P}{{\bf P}}
\begin{document}
\centerline{\bf
 ON KOSZUL PROPERTY OF THE HOMOGENEOUS}
\centerline{\bf  COORDINATE RING OF A CURVE}
\bigskip
\centerline{Alexander Polishchuk}
\bigskip
\bigskip

\section{Introduction}

This paper is devoted to Koszul property of the homogeneous
coordinate algebra of a smooth complex algebraic curve in the
projective space (the notion of a Koszul algebra
 is some homological refinement of the notion of a quadratic algebra,
for precise definition see next section). It
grew out from the attempt to understand methods of M. Finkelberg and
A. Vishik in their paper \cite{FV} proving this property for the
canonical algebra of a curve in the case it is quadratic. The basic
ingredient of their proof is the following lemma on special divisors.

\begin{lem}{\rm (\cite{GL})} Let $C$ be a non-hyperelliptic
non-trigonal curve which is not a plane quintic. Then
 there exists a divisor $D$ of degree $g-1$ on $C$
such that $|D|$ and $|K(-D)|$ are base-point free linear systems of
dimension 1 where $K$ is the canonical class.
\end{lem}

In the cited paper of M. Green and R. Lazarsfeld it is used for
 the vector bundle proof of Petri's theorem which asserts that if $C$
satisfies the conditions of this lemma then the canonical algebra of $C$
is quadratic.
 In this note we will show that Koszul property can be
derived from this lemma by purely homological technique combined with
a simple statement concerning Koszul property of the homogeneous
coordinate algebra of a finite set of points in the projective space
--- in particular we obtain a new proof of Petri's theorem. It turns
out that the same technique works for the proof of Koszul property of
the homogeneous coordinate ring of a curve of genus $g$ embedded by
complete linear system of degree $\ge 2g+2$ (this result is due to D.
Butler \cite{B}) and also for embeddings defined by the
complement to some very special linear systems (the simplest case being
that of tetragonal systems) in the canonical class.
In particular, we prove Koszul property for a general tetragonal
curve of genus $g\ge 9$ embedded by $K(-T)$ where $T$ is the tetragonal
series. It is worthy to mention here
that if $C$ is a non-hyperelliptic curve and $L$ is a linear bundle
of degree $2g+1$ such that $H^0(L\ot K^{-1})=0$ then the coordinate algebra of
$C$ in the
embedding defined by $L$ is quadratic (see e.g.
\cite{L}). It is still
unknown (at least to me) whether it is Koszul or not. It seems that
the following general question is also open: whether the homogeneous
coordinate algebra of a projectively normal smooth connected
complex curve is Koszul provided it is quadratic?

To illustrate our technique we present here the proof of Petri's theorem
based on the lemma mentioned above.

\begin{thm} Assume that there exists a divisor $D$ on $C$ of degree
$g-1$ such that $|D|$ and $|K(-D)|$ are base-point-free linear series
of dimension 1.  Then the canonical algebra $R$ is quadratic.
\end{thm}

\Pf . Choose the divisor $P_1+\ldots+P_{g-1}$ in the linear system $|D|$
such that these $g-1$ points are distinct (then they are
in general linear position in the
$(g-3)$-dimensional space they span). Let $V=H^0(K(-D))\subset
H^0(K)=R_1$ be the corresponding 2-dimensional subspace. Denote the
symmetric algebra $Sym(R_1)$ by $S$ and let $J\subset S$ be the
homogeneous ideal of $C$ so that $R=S/J$ by Noether's theorem
(obviously $C$ is non-hyperelliptic ).  Consider the following exact
sequence $$0\ra VS\cap J\ra J\ra J/(VS\cap J)\ra 0$$ where $VS$ is the
ideal in $S$ generated by $V$.  It is easy to see that it is enough to
check the following two statements:
\begin{enumerate}
\item
$J/(VS\cap J)$ is generated
over $S$ by elements of degree 2;
\item
$VS\cap J$ is generated over $S$
by elements of degree 2 modulo $VJ$.
\end{enumerate}
We will use the following lemma.

\begin{lem} Let $|D|$ be a base-point-free linear system of dimension 1.
 Then the natural homomorphism $H^0(D)\ot H^0(K^n)\ra H^0(K^n(D))$
is surjective for $n\ge 0$.
\end{lem}

The proof is left to reader.

Now for the proof of 1) we claim that $VS+J$ is exactly the
homogeneous ideal of the set of points $P_1,\ldots, P_{g-1}$ so the
statement follows from the quadratic property of an algebra of $g-1$
points in general position in $\P^{g-3}$ (see sect. 3 for more general
result in this direction). Indeed let $A$ be the coordinate algebra of
these points. Then $A$ is the factoralgebra of $R$ by the ideal
$I=\oplus_{n\ge1}H^0(K^n(-D))$. Now it follows from the above lemma
that $I$ is generated by $I_1=V$ over $R$ so $A=R/(VR)=S/(VS+J)$ and
our assertion follows.

It remains to check 2). For this consider the multiplication map
\\ $\mu_S:V\ot S(-1)\ra S$. It is enough to verify that
$\mu_S^{-1}(J)/(V\ot J(-1))$ is generated over $S$ by elements of degree 2. We
can
rewrite this as follows: let $\mu_R:V\ot R(-1)\ra R$ be another multiplication
map then ker$\mu_R$ is generated by elements of degree 2 over $R$. But this
kernel is equal to $\oplus_{n\ge2}H^0(K^{n-2}(D))$ (by base-point-pencil trick
applied to $K(-D)$) so we are done by the lemma above. \qed

{\it Acknowledgments.} I am grateful to A. Bondal, L. Positselsky and
A. Vishik for stimulating discussions. Also I thank J. Kollar and
University of Utah for their hospitality while carrying out this
research.

\section{Some homological algebra}

In this section we prove a general criterion for Koszul property of a
graded algebra given its Koszul factoralgebra and some information
about this factoralgabra as a module over original algebra.

Recall that graded (associative but not necessary commutative) algebra
$A=A_0\oplus A_1\oplus\ldots$ over
the field $k$ is called Koszul if $A_0=k$ and ${\rm Ext}^n(k,k(-m))=0$
for $n\neq m$. Here functor Ext is taken in the category of graded
(left) $A$-modules, $k$ is trivial $A$-module concentrated in degree
0, for graded $A$-module $M=\oplus M_i$ we define the shifted module
as $M(l)=\oplus M_{i+l}$. Note that the above condition on Ext's
for $n=1$ means that $A$ is generated by elements of degree 1,
for $n=1$ and $n=2$ - that $A$ is quadratic (i.e. in addition it has
defining relations of degree 2). For equivalent definitions via
exactness of Koszul complex and distributivity of some lattices see
\cite{P},\cite{Ba1},\cite{BGS}.

The following theorem is the generalization of lemma 7.5
 in \cite{BP} where
it was proved in the particular case when $A=k$.

\begin{thm} Let $R$ be graded algebra over $k$ with $R_0=k$ and $A$ --
its Koszul (graded) factoralgebra. Assume that there exists a complex
$K^{\cdot}$ of free right $R$-modules of the form
$$\ldots\ra V_2\ot R(-2)\ra V_1\ot R(-1)\ra R$$
(so $K^i=V^i\ot R(-i)$ where $V^i$ are finite-dimensional
$k$-linear spaces, $V^0=k$) such that $H_0(K^.)=A$, $H_p(K^.)_j=0$ for
$p\ge 1$, $j>p+1$. Then $R$ is Koszul and $A$ has a linear free
resolution as $R$-module that is a resolution of the form
$$\ldots\ra U_2\ot R(-2)\ra U_1\ot R(-1)\ra R\ra A\ra 0$$
\end{thm}

\Pf . Note that the conditions on the homology $H_p(K^{\cdot})$
for $p\ge 1$ in formulation of the theorem mean that $H_p(K^{\cdot})$ as
$R$-module is an extension of some multiples of trivial modules
$k(-p)$ and $k(-p-1)$. Indeed this follows directly from the fact that
$H_p(K^{\cdot})_j=0$ for $j<p$ which is evident. Now we prove by induction
on $n\ge 0$ that $Ext_R^n(k,k(-m))=0$ if $m\neq n$. Case $n=0$ is
trivial. Assume that $n>0$, $\E_R^i(k,k(-m))=0$ for $i<n$, $m\neq
i$ and let us prove the assertion for $n$. Considering the spectral
sequence
$$E_2^{p,q}=\E_A^q(Tor_p^R(A,k),k(-m))\Ra\E_R^{p+q}(k,k(-m))$$ it is
easy to see that it is enough to prove that ${\rm Tor}_j^R(A,k)_i=0$
if $i\neq j$, $j\le n$ (here we use Koszul property of $A$).  The
latter is equivalent to $\E_R^j(A,k(-i))=0$ for $i\neq j$, $j\le n$
(here $\E$ is taken in the category of right $R$-modules). Consider
the spectral sequence associated with complex $K^{\cdot}$ and
cohomological functor ${\rm Hom}_{D^-(R)}(.,k(-m))$ where $D^-(R)$ is
the derived category of complexes bounded from the right:
$$E_2^{p,q}=\E_R^q(H_p(K^{\cdot}),k(-m))\Ra
{\rm Hom}^{p+q}_{D^-(R)}(K^{\cdot},k(-m))$$
The limit on the right can be computed easily with help of another
spectral sequence with $E_1$ obtained by applying the same
functor to the terms of the complex $K^{\cdot}$. Due to the form of
this complex it degenerates at $E_1$ which gives an equality ${\rm
Hom}^n(K^{\cdot},k(-m))=0$ if $n\neq m$.  On the other side for
$p\ge1$ the group $H_p(K^{\cdot})$ is an extension of the direct sums
of several copies of $k(-p)$ and $k(-p-1)$ and we obtain by assumption
that $E_2^{p,q}=\E_R^q(H_p(K^{\cdot}),k(-m))=0$ if $p\ge1$, $q<n$,
$m\neq p+q$, $m\neq p+q+1$. This implies that the terms $E_2^{0,q}$
must survive and contribute to $E_{\infty}$ if $q\le n$, $q\le m-1$.
Indeed all the differentials $d_r: E_r^{r-1,q-r}\ra E_r^{0,q}$ are
zero for $r\ge2$, $q\le {\rm min}(n,m-1)$ because $E_r^{r-1,q-r}$ is
zero for these values of $r$ and $q$. Therefore our computation of the
limit implies an equality $E_2^{0,q}=\E_R^q(A,k(-m))=0$ for $q\le{\rm
min}(n,m-1)$. On the other side obviously $\E_R^q(A,k(-m))=0$ for
$q>m$ so we are done. \qed

\section{Application to the coordinate ring of a curve}

Now we are going to apply the theorem of the previous section to the
homogeneous coordinate algebra of a curve. So let $C$ be a curve, $L$
be a very ample linear bundle on $C$. We are interested in the algebra
$R=R_L=\\ =\oplus_{n\ge0}H^0(L^n)$. Our approach is as follows: we
consider the factoralgebra $A$ of $R$ associated with some effective
divisor $D$ on $C$ namely $A=A_D=R/J_D$ where
$J_D=\oplus_{n\ge1}H^0(L^n(-D))$ is an ideal in $R$. To apply the
above theorem to this situation we have to construct a complex of free
$R$-modules which is an "almost resolution" of $A$ and has the
required form, and to check Koszul property of the algebra $A$. The
latter is simple provided that points of the divisor $D$ are
"sufficiently linear independent" in the embedding defined by $L$. To
construct the desired complex we assume that $\O(D)$ and $L(-D)$ are
base-point-free so that we have the following exact triples:  $$0\ra
 \O(-D)\ra V\ot \O\ra \O(D)\ra0$$ $$0\ra L^{-1}(D)\ra U\ot\O\ra
L(-D)\ra0$$ where $V$ and $U$ are some vector space of dimension 2.
Now the complex $K$ has the following form:
$$\ldots \ra U\ot R(-3)\stackrel{d_3}{\ra} V\ot R(-2)\stackrel{d_2}{\ra} U\ot
R(-1)\stackrel{d_1}{\ra} R$$
where $d_1$ is induced by the composition of maps
$$U\ra H^0(L(-D))\ra H^0(L)=R_1,$$
the $n$-th component of $d_{2k}$ is the composition
$$V\ot H^0(L^{n-2k})\ra H^0(L^{n-2k}(D))\ra U\ot H^0(L^{n-2k+1}),$$
and that of $d_{2k+1}$ is the composition
$$U\ot H^0(L^{n-2k-1})\ra H^0(L^{n-2k}(-D))\ra V\ot H^0(L^{n-2k})$$
Using the exact triples
above one can compute easily the homology of $K$. Indeed
$${\rm ker}(d_{2k+1})_n={\rm ker}(U\ot H^0(L^{n-2k-1})\ra
H^0(L^{n-2k}(-D)))\simeq H^0(L^{n-2k-2}(D))$$
Therefore
$$H_{2k+1}(K)_n\simeq{\rm coker}(V\ot H^0(L^{n-2k-2})\ra
H^0(L^{n-2k-2}(D)))$$ $$\simeq{\rm ker}(H^1(L^{n-2k-2}(-D))\ra V\ot
H^1(L^{n-2k-2})$$
It follows that if $H^1(L^2(-D))=0$ then
$H_{2k+1}(K)_n=0$ for $n\ge 2k+4$. Futhermore if the map
$\alpha:H^1(L(-D))\ra V\ot H^1(L)$ is injective then
$H_{2k+1}(K)_{2k+3}=0$ as well. In analogous way we obtain that if
$H^1(L(D))=0$ and $\beta:H^1(D)\ra U\ot H^1(L)$ is injective then
$H_{2k}(K)_{\ge 2k+2}=0$. Note that the condition $H^1(L(D))=0$
implies surjectivity of $\alpha$ so in this case injectivity is
equivalent to the equality of dimensions: $h^1(L(-D))=2h^1(L)$.
Analogously if $H^1(L^2(-D))=0$ then injectivity of $\beta$ is
equivalent to equality $h^1(D)=2h^1(L)$. It easy to see that if all
these conditions hold and in addition $h^0(L(-D))=2$ that is $U\simeq
H^0(L(-D))$, then image of $d_1$ is equal to the ideal $J_D$; thus all
the homological conditions of Theorem 1 are satisfied except for
Koszul property of an algebra $A$. At this step we use the following
result of G. Kempf.

\begin{thm}{\rm (\cite{K})}
 Homogeneous coordinate algebra of the finite set of
$d$ distinct points in general linear position in $\P^{d-p}$
is Koszul provided that $p\le d/2$.
\end{thm}

\begin{rem} In the case $p\le 3$ this is particularly easy: the points lie
on a rational normal curve (see \cite{GH}) so the statement can
be deduced easily from the main theorem of \cite{Ba2} (see also \cite{ERT}).
Note also that if $p=1$ (resp. $p=2$) then the set of points in
question is a hyperplane section of a rational (resp. an elliptic)
normal curve so the Koszul property in these two cases
follows from the same property of the coordinate algebra of a rational
(resp. an elliptic) normal curve.  So we can avoid the reference to
the results above if we are interested in the case $h^1(L)\le 1$ only.
\end{rem}

Now we are ready to prove our main theorem.

\begin{thm} Let $L$ be a very ample linear bundle on $C$ of degree
\\ ${\rm deg}L\ge g+3$ such that corresponding embedding of $C$ into
$\P(H^0(L)^*)$ is projectively normal. Assume that
there exists a divisor $D=P_1+\ldots+P_d$ of the degree
$d={\rm deg}L-g-1+2h^1(L)$
such that the linear series $|D|$ and $|L(-D)|$ are base-point free of
the dimensions ${\rm deg}L-2g+4h^1(L)-1$ and 1 correspondingly. Assume
also that $h^1(L(D))=h^1(L^2(-D))=0$ and that
any $h^1(L)$ points of $D$ impose
independent conditions on $K\ot L^{-1}(D)$.
Then algebra $R_L$ is Koszul.
\end{thm}

\Pf . By Riemann-Roch we obtain that $h^1(D)=h^1(L(-D))=2h^1(L)$ so
it follows from the discussion
above that we only have to check Koszul property of the set of points
${P_1,\ldots,P_d}$ embedded by $L$. The
dimension of the linear subspace spanned by these points is equal to
$h^0(L)-3=d-h^1(L)-1$. The condition $p=h^1(L)+1\le d/2$ is satisfied by
 assumption
so by Theorem 5 it is sufficient to
verify  that any
$d-h^1(L)$ of $P_1, \ldots, P_d$ impose independent conditions on
$|L|$. But this is equivalent to the property we have assumed that any
$h^1(L)$ of them impose independent conditions on $K\ot L^{-1}(D)$.
\qed

\begin{rem} The condition of the theorem concerning independence of
any $h^1(L)$  points  is  satisfied  automatically if $h^1(L)\le 1$.
If $h^1(L)=2$ and ${\rm deg}L=2g-6$ so that the dimension of  $|D|$
is  1 then  the  sufficient  condition  is that $K\ot L^{-1}(2D)$ is
very ample.
\end{rem}

\begin{cor} If $C$ is a non-hyperelliptic, non-trigonal curve which
is not a plane quintic then the canonical algebra $R_K$ is Koszul.
\end{cor}

\Pf . It follows from the Green-Lazarsfeld's lemma and the theorem above. \qed

\begin{cor} For any curve of genus $g$ and any linear bundle $L$ of
degree $\ge 2g+2$ algebra $R_L$ is Koszul.
\end{cor}

\Pf . In this case $h^1(L)=0$ and we can choose $D$ as above
in the linear system $L(-P_1-\ldots -P_{g+1})$ for general $g+1$
points $P_1,\ldots ,P_{g+1}$. Note also that Koszul property in this
case follows trivially from Kempf's Theorem applied to a general
hyperplane section of $C$. \qed

\begin{cor} Under the assumptions of the theorem the following natural map
is surjective:
$$H^0(L)\ot H^0(D)\ra H^0(L(D))$$
Also if we define a vector bundle $M_D$ from the exact triple
$$0\ra M_D\ra H^0(D)\ot \O\ra \O(D)\ra 0$$
then $M_D\ot L$ is generated by global sections.
\end{cor}

\Pf . Both statements follow from the second part of Theorem 4. \qed

It remains to analyze the case $h^1(L)\ge 2$ of our theorem.
Put $L=K(-A)$, dim$|A|=r$.
Then the condition $h^0(D)\ge 2$ implies the following
inequality: deg$A\le 4r$. On the other hand
considering the natural map
$$|A|\times |D|\ra |K\ot L^{-1}(D)|$$
we obtain another inequality:
$$r+h^0(D)-1\leq h^1(L(-D))-1$$
which is equivalent to deg$A\ge 3r$.
Also considering the map
$$|D|\times|L(-D)|\ra |L|$$
we obtain the restriction $r\le g/3-1$.
In the case $r=1$ we obtain from the above inequalities that $|A|$ is
either trigonal or tetragonal system. Furthermore it is easy to see that
the former case is impossible so $A=T$ is  a tetragonal system.
One can check that $T$ should be base-point-free otherwise $C$
fails to be cut out (even set-theoretical) by quadrics in the
embedding defined by $|L|$. Also the fact that $L$ is very ample implies
that $C$ is not hyperelliptic.
To satisfy
the conditions of the theorem we should have a decomposition of $L$
into the sum of two divisors of degree $g-3$ each defining the
base-point free linear system of dimension 1. In the next section we
will give some examples when this situation occurs.

\begin{rem} So far I don't know much about the case $h^0(L)>2$.
I hope that there should be examples when the above technique
applies to this case too.  \end{rem}

\section{Tetragonal curves}

In this section we study the case of the embedding of a tetragonal
curve $C$ by the complete linear system $|K(-T)|$ where $T$ is a
(base-point-free)
tetragonal series. As a tetragonal series is not unique in general it
is natural to consider the moduli space ${\cal M}_g^t$ of pairs $(C,
T)$ where $T$ is such a series on a non-hyperelliptic curve $C$ of
genus $g$.  By the well-known construction (see \cite{S}) $T$ gives an
embedding of $C$ into 3-dimensional rational normal scroll $X$.
Furthermore it is easy to see that $C$ is a complete intersection of
two divisors on $X$. Thus there is a stratification of ${\cal M}_g^t$
by the type of scroll $X$ and the type of complete
intersection, and all the strata are irreducible
{}.
Note that the Hilbert series of $R_{K(-T)}$ is constant over ${\cal M}_g^t$
so by the well-known result the set of pairs $(C,T)$ for which $R_{K(-T)}$ is
Koszul is an intersection of countably many open subsets in ${\cal M}_g^t$.
The problem to be solved is to find all the
strata ${\cal N}$ such that for general pair $(C,T)\in {\cal N}$
the algebra $R_{K(-T)}$ is Koszul (we say that such a stratum is Koszul).
Here
"general" means "in the complement of countably many proper
subvarieties" so a stratum satisfies this property if it contains at
least one such pair. Note that if a stratum ${\cal N}_1$ is Koszul and
it is contained in the closure of a stratum ${\cal N}_2$ then ${\cal N}_2$
is also Koszul. The analogous problem for quadratic algebras is easy and
we will see that the most degenerate (with respect
to complete intersection type) quadratic strata are Koszul,
so it is natural to expect that this is true in general.
Also we construct examples of pairs $(C,T)$, for which
algebra $R_{K(-T)}$ is Koszul, on some other strata
cosidering ramified double coverings of
hyperelliptic curves and applying the method of the previous section.
As an evident consequence we
obtain that for general tetragonal curve of genus $g\ge 9$
algebra $R_{K(-T)}$ is Koszul
(note that for general tetragonal curve $T$ is unique).

We begin with recalling the construction of a 3-dimensional rational
normal scroll containing tetragonal curve $C$. Let $T$ be a
base-point free tetragonal system on $C$. It defines a 4-sheeted
covering $\pi:C\ra\P^1$. Applying the relative duality to the
canonical non-vanishing section $\O\ra \pi_*\O$ we obtain the
surjective homomorphism $\pi_*K_C\ra\O(-2)$. Let $V$ be its kernel so
that we have the following exact triple on $\P^1$:
$$0\ra V\ra \pi_*K_C\ra \O(-2)\ra 0$$
It follows that the homomorphism $V\ra
\pi_*K_C$ induces an isomorphism of global sections and therefore
the corresponding homomorphism $\pi^*V\ra K_C$ is surjective. Thus we
obtain a morphism $\phi:C\ra X$ where $X=\P(V^{\vee})$ such that
$\phi^*\O_X(1)\simeq K_C$.  Note that it follows from the exact triple
 above that $h^0(V(-i))=h^0(K_C(-iT))$. In particular $h^0(V)=g$,
$h^0(V(-1))=g-3$ so that $V\ge0$ in the  sense that all linear direct
summands in $V$ has form $\O(l)$ with $l\ge0$. Therefore $\O_X(1)$ is
base-point free and defines a morphism from $X$ to $\P^{g-1}$ inducing
the canonical morphism by composition with $\phi$. It follows that
$\phi$ is an embedding (we have assumed that $C$ is
non-hyperelliptic).  Furthermore $h^0(V(-2))=g-6$ if and only if
$h^0(2T)=3$ and assuming that we obtain that $V(-1)\ge0$ and hence
$\O_X(1)$ is very ample.

 Let
$p:X\ra \P^1$ be the projection.  The push forward by $p$ of the
natural homomorphism $\O(2)\ra \O(2)|_C$ gives rise to a homomorphism
$f:S^2V\ra\\ \ra \pi_*K_C^2$ where $S^2V\simeq p_*\O(2)$ is the second
symmetric power of $V$. We claim that $f$ is surjective. Indeed it
suffices to check this pointwise so we have to verify that for each
$q\in \P^1$ the four points of $\pi^{-1}(q)$ impose independent
conditions on quadrics in the corresponding projective plane
$p^{-1}(q)$. But this is evidently true because these points span that
plane by the geometric form of the Riemann-Roch Theorem. Let us denote
the kernel of $f$ by $E$. This is a rank-2 vector bundle on $\P^1$ of
degree $g-5$ and it fits in the following exact sequence
$$0\ra E\ra S^2V\ra \pi_*K^2\ra0$$
Now we claim that the corresponding homomorphism $p^*E\ra
\O_X(2)$ vanishes exactly along $C\in X$. Indeed as any divisor in
$|T|$ contains four points no three of which lie on a line they are
cut out by two quadrics in the corresponding plane. More than that,
comparing arithmetical genera we
conclude that $C$ as a subscheme of $X$ coincides with the zero-locus
of the corresponding regular section of $p^*E^{\vee}(2)$. Now $E\simeq
\O(a)\oplus\O(b)$ where $a+b=g-5$ so $C$ is in fact a complete
intersection of two divisors $S_1\in|\O_X(2)(-aH)|$ and
$S_2\in|\O_X(2)(-bH)|$ where $H=p^*\O(1)$.

Assume now that $K_C(-T)$ is projectively normal. Then the
homomorphism $S^2H^0(K_C(-T))\ra H^0(K_C^2)$ is surjective. As
$V(-1)\ge0$ it follows that the natural homomorphism $S^2H^0(V(-1))\ra
H^0(S^2V(-2))$ is surjective too. Hence the homomorphism $S^2V(-2)\ra
\pi_*(K_C^2)(-2)$ induces a surjection on global sections and
consequently $H^1(E(-2))=0$ that is $a,b\ge1$. Conversely it easy to
see that if $C$ is a complete intersection as above with $a,b\ge1$
then $K_C(-T)=(\O_X(1)(-H))|_C$ is projectively normal. We fix these
results in the following proposition.

\begin{prop} Let $C$ be a non-hyperelliptic curve of genus $g$ with a
base-point free tetragonal system $T$ on it. Then we can present
$C$ as a complete intersections of two divisors from the linear
systems $|\O_X(2)(-aH))|$ and $|\O_X(2)(-bH)|$  on a 3-dimensional
rational normal scroll $X=\P(V^{\vee})$, where $V\ge 0$ is rk-3 vector
bundle of degree $g-3$ on $\P^1$, in such a way that $T=H|_C$ (here
$H$ is the pull-back of $\O(1)$ from $\P^1$, $a+b=g-5$). Furthermore
$h^0(2T)=3$ if and only if $V(-1)\ge0$ and if this condition is
satisfied then $K_C(-T)$ is projectively normal if and only if $a\ge1,
b\ge1$.
\end{prop}

\begin{rem} Note that if $K_C(-T)$ is very ample then $a,b\ge0$.
In any case $a,b\ge -1$.
\end{rem}

Our next remark is that under notations of the previous proposition
$R_{K_C(-T)}$ is quadratic if and only if $a,b\ge 2$ (this is very easy
to verify using exact sequences of the restriction to a divisor). So
it makes reasonable the following conjecture.

\begin{conj} Under assumptions and notations above
if  $a,b\ge 2$ then the algebra $R_{K_C(-T)}$ is Koszul .
\end{conj}

\begin{rem} It is easy to see that this is true at least when $a=2$ or $b=2$.
Indeed first we can prove that the coordinate algebra of any divisor
$S\in |\O(2)(-aH)|$ under the embedding defined by $\O(1)(-H)$ is Koszul
provided that $a\ge 2$. The reason is that its hyperplane section is
a curve of genus $g-5-a$ embedded by the complete linear system of degree
$2g-10-a$ so we can apply the corollary 2 above. Then our curve $C$ is an
intersection of $S$ with a quadric so its coordinate algebra is Koszul
by the result of \cite{BF}. This proves the conjecture for
$g=9,10$. To prove the Koszul property for general pair $(C,T)$ of some
stratum with $a,b\ge 2$
it would be sufficient
to prove that this stratum can be degenerated into one with $a=2$ or $b=2$.
\end{rem}

Our method of proving Koszul property suggests the following conjecture
which implies the previous one.

\begin{conj} Under the same assumptions there exists a decomposition
of $K_C(-T)$ into the sum of two base-point free pencils of
degree $g-3$.
\end{conj}

Now we consider the specific case when our tetragonal curve $C$ is a
double covering of a hyperelliptic curve. We are going to construct
some examples when the required decomposition of $K_C(-T)$ into the
sum of  two pencils exists. For this we will use the well-known
connection between linear bundles over the double covering and rk-2
bundles with a Higgs field (which is a twisted endomorphism of a
bundle) over the base of the covering (see for example \cite{Hi}).

So let $C_h$ be the hyperelliptic curve of genus $g_h$,
$\pi_h:C_h\ra\P^1$ be the corresponding double covering,
$\G=\pi_h^*\O(1)$ be the hyperelliptic system. Then
$K_h=\O((g_h-1)\G)$ is the canonical class of $C_h$.
Now we consider a Higgs bundle $(F,\phi)$ where $F=\O(D_1)\oplus
\O(D_2)$, $D_i$ being the divisors of degree $g_h-2$ and $\phi:F\ra
F(M)$ is a homomorphism defined by the sections $s_1\in
H^0(M(D_2-D_1))$ and $s_2\in H^0(M(D_1-D_2))$ (other enties of
$\phi$ being zero) -- here M is some linear bundle which has
sufficiently large degree to be estimate later. These data define the
line bundle $\O(D)$ over the double covering $\pi:C\ra C_h$ such that
$\pi_*\O(D)\simeq F$. Simple computation shows that deg$D=g-3$ where
$g=2g_h-1+{\rm deg}M$ is the genus of $C$. Now we look under what
conditions $|D|$ is  a base-point-free linear system of dimension 1.
First we should have that $F$ is globally generated outside the
ramification divisor $s=s_1s_2\in H^0(M^2)$. Hence
$h^0(D_1)=h^0(D_2)=1$ and the unique divisor of $|D_i|$ is contained in
the zero divizor of $s$. Furthermore at the point $x$ of ramification
global sections of $F$ should generate the unique $\phi$-invariant
1-dimensional factorspace of the stalk $F_x$. It is easy to see that
these conditions are satisfied if we put $s_i=u_it_i$ where $u_i\in
H^0(D_i)$ are non-zero sections, $(i=1,2)$, $t_1\in H^0(M(D_2-2D_1))$,
$t_2\in H^0(M(D_1-2D_2))$ provided that zero divisors of $u_1, u_2,
t_1, t_2$ are all disjoint. Note that the change of $D$ by $K_C(-T-D)$
where $T=\pi^*\G$ leads to the change of $D_i$ by $(g_h-2)\G-D_i$
with essentially the same Higgs field. In particular if we choose
$D_2=(g_h-2)\G-D_1$ then we will have $\O(2D)\simeq K_C(-T)$. Now it is
clear that if $M$ is a general linear bundle of degree at least
$2g_h-1$ then we can find $(F,\phi)$ as above such that corresponding
divisor $D$ on $C$ satisfies the conditions of Theorem 2. Indeed
then deg$M(D_2-2D_1)\ge g_h+1$ and we can use the fact that general
bundle of degree $\ge g_h+1$ is base-point free. Also as we have
mentioned above under suitable choices we'll have $\O(2D+T)\simeq K_C$
which implies the last condition of Theorem 2 (see remark after it).

Now in order to verify projective normality of $K_C(-T)$ we have to
compute the discrete invariants described above (namely the bundles
$V$ and $E$ on $\P^1$) of the pair $(C,T)$ obtained in this way. First
we have the canonically splitting exact triple
 $$0\ra K_h\ot M\ra \pi_*K_C\ra K_h\ra 0$$
 Now $V$ is the kernel of composition
$$(\pi_h)_*(\pi_*K_C)\ra(\pi_h)_*K_h\ra \O(-2)$$
 It follows that $V$ fits into the splitting exact sequence
$$0\ra (\pi_h)_*M(g_h-1)\ra V\ra \O(g_h-1)\ra 0$$
Now it is easy to check that $(\pi_h)_*\pi_*(K_C^2)\simeq((\pi_h)_*(M^2)\oplus
(\pi_h)_*M)(2g_h-2)$ and the natural map
$$ S^2((\pi_h)_*M\oplus \O)\simeq S^2V(-2g_h+2)\ra
(\pi_h)_*\pi_*(K_C^2)(-2g_h+2)\simeq$$ $$ \simeq
(\pi_h)_*(M^2)\oplus(\pi_h)_*M$$
is induced by the natural maps $\psi:S^2((\pi_h)_*M)\ra (\pi_h)_*(M^2)$,
$\O\ra (\pi_h)_*(M^2)$ and the identity map of $(\pi_h)_*M$. It follows that
there is an exact sequence
$$0\ra {\rm ker}\psi \ra E(-2g_h+2)\ra \O\ra 0$$
provided that $\psi$ is surjective. Assume that deg$M\ge 2g_h+1$ then $\psi$ is
evidently surjective and ker$\psi\simeq\O({\rm deg}M-2g_h-2)$ hence we obtain
that $E(-2g_h+2)\simeq \O\oplus \O({\rm deg}M-2g_h-2)$. At last note that as
$M$ is general the splitting type of $(\pi_h)_* M$ is either $(i,i)$
or $(i,i+1)$ depending on the parity of deg$M$.

	 Summarizing the discussion above we obtain the
following statement.

\begin{thm} Let $\cal N$ be the stratum of the moduli space ${\cal
M}_g^t$ $(g\ge 9)$ with one of the following splitting types of $V$
and $E$:  \begin{enumerate} \item $(g_h-1,g_h-1+i,g_h-1+i)$ and
$(2g_h-2,g_h-3+2i)$ where $g=3g_h+2i$, $i\ge g_h/2$, $g_h\ge2$;
\item  $(g_h-1,g_h-2+i,g_h-1+i)$ and $(2g_h-2,g_h-4+2i)$ where
$g=3g_h+2i-1$, $i\ge(g_h+1)/2$, $g_h\ge2$.  \end{enumerate}
Then for general $(C,T)\in{\cal N}$ algebra $R_{K(-T)}$ is Koszul.
\end{thm}

\begin{rem} The condition $g\ge9$ excludes the case $g_h=2, i=1$ in
1).  \end{rem}

\begin{cor}
For general tetragonal curve of genus $g\ge 9$ algebra $R_{K(-T)}$ is Koszul.
\end{cor}

\Pf . It is sufficient to note that the pairs $(C,T)$ constructed above
(for which this algebra is Koszul) lie in the open subset $U$
of ${\cal M}_g^t$ for which $a,b\ge 1$ (notation as above), and it is known
(see \cite{S}) that for such curves $C$ the series $T$ is unique. Hence
$U$ embeds as an open subset in the locus of tetragonal curves inside
${\cal M}_g$ which is irreducible, so $U$ itself is irreducible.
\qed

\section{On regularity of modules over a commutative Koszul algebra}

In this section we give a geometric bound for the regularity of a
module over a commutative Koszul algebra.
The result is not new but our proof seems to be very simple so we
present it here.

Let $X$ be a projective scheme, $L=\O_X(1)$ be a very ample line bundle
on $X$ such that the corresponding algebra $R=R_L$ is Koszul. For a sheaf
$F$ on $X$ we denote tensor product of $F$ with the $n$-th power of $L$
by $F(n)$. With these assumptions we have the following.

\begin{thm} For any coherent sheaf $F$ on $X$ if $H^i(F(-i))=0$ for any $i>0$
then the corresponding $R$-module $M=\oplus_{i\ge 0}H^0(F(i))$ has a
linear free resolution that is a resolution of the form
$$\ldots \ra V_2 \ot R(-2) \ra V_1 \ot R(-1) \ra V_0 \ot R \ra M \ra 0$$
where $V_i$ are some vector spaces.
\end{thm}

\Pf . Koszul property of $R$ means that there is a resolution of the trivial
module (of degree zero) which has form
$$\ldots \ra Q_2 \ot R(-2) \ra Q_1 \ot R(-1) \ra Q_0\ot R \ra k \ra
0$$ where $Q_0=k$.
It induces the following exact sequence of sheaves on $X$ $$\ldots
\ra Q_2 \ot \O(-2) \ra Q_1 \ot \O(-1) \ra \O \ra 0$$ Now we can tensor
it by $F(n)$ and consider the corresponding spectral sequence
computing hypercohomology
$$E_1^{p,q}=Q_{-p}\ot H^q(F(i+p))\Rightarrow 0$$
with differentials $d_r$ of bidegree $(r,-r+1)$.
Now $F$ is $0$-regular in the sense of Castelnuovo-Mumford (see
\cite{M}) so we have $E_1^{p,q}=0$ if $p+q\ge -i,\ q\ge 1$.
It follows that the complex $E_1^{\cdot,0}$ is exact in terms $p>-i$.
But this complex computes syzigies of $M$, namely its cohomology in
$p$-th term is Tor$^R_p(k,M)_(i+2p)$. So we have proved that
Tor$_p(k,M)_j=0$ for $j>p$ which is equivalent to the existence of
linear free resolution in question. \qed

\begin{rem} The statement of this theorem is
essentially equivalent to the statement of the main theorem of
\cite{AE} for the case of Koszul algebras $R$ which has form $R_L$.
Indeed that theorem asserts that the regularity of a module $M$
over $R$ is bounded by its regularity as a module over a symmetric
algebra $S$ which surjects onto $R$. We have proved that for a modules which
come from coherent sheaves the regularity is bounded by its geometric
counterpart --- the regularity in the sense of Castelnuovo-Mumford.
 Note that the regularity of an arbitrary module (not
necessary coming from coherent sheaf) can be bounded easily using the
bound for that special type of modules.
Now the point is that for the case of symmetric algebra these two
regularities are equal so we arrive to the formulation of \cite{AE}.
\end{rem}

%


Harvard University,

{\it e-mail}: apolish@math.harvard.edu

\end{document}